\documentclass{myarticle}

\usepackage{txfonts}
\usepackage[colorlinks]{hyperref}
\usepackage{tikz}

\usetikzlibrary{positioning}

\renewcommand{\d}{\displaystyle}
\newcommand{\n}{\noindent}
\newtheorem{thm}{Theorem}
\newtheorem{prop}{Proposition}

\newtheorem{cor}{Corollary}

\begin{document}

\title{Connective eccentric index of some graph operations}

\author[nd]{Nilanjan De}
\ead{de.nilanjan@rediffmail.com}
\author[ap]{Anita Pal}
\ead{anita.buie@gmail.com}
\author[man]{Sk. Md. Abu Nayeem\corref{cor1}}
\ead{nayeem.math@aliah.ac.in}

\address[nd]{Department of
Basic Sciences and Humanities (Mathematics),\\ Calcutta Institute of Engineering and Management, Kolkata, India.}
\address[ap]{Department of
Mathematics, National Institute of Technology, Durgapur, India.}
\address[man]{Department of Mathematics, Aliah University,\\ DN - 20, Sector -  V, Salt Lake, Kolkata - 700 091, India.}
\cortext[cor1]{Corresponding Author.}

\begin{abstract}
The connective eccentric index of a graph is a topological index involving degrees and eccentricities of vertices of the graph. In this paper, we have studied the connective eccentric index for double graph and double cover. Also we give the connective eccentric
index for some graph operations such as joins, symmetric difference, disjunction and
splice of graphs.

\end{abstract}

\maketitle

\section{Introduction}
Let $G$ be a simple connected graph with vertex set $V(G)$ and edge set $E(G)$. Let
$n$ and $m$ be the number of vertices and edges of $G$
respectively. We denote the degree of a vertex $v$ of $G$ by $d_G(v)$. For $u,v \in V(G)$, the distance between $u$
and $v$ is defined as the length of a shortest path between $u$
and $v$. For a given vertex $v$ of $G$, the eccentricity $\varepsilon_G(v)$ is the largest distance from $v$ to any other vertices  of $G$. The sum of eccentricities of all the
vertices of $G$ is denoted by $\theta(G)$ \cite{dan04}. If any vertex $v \in V(G)$ is
adjacent to all the other vertices of $G$ then $v$ is called a
well-connected vertex. Thus, if $v \in V(G)$ is a well-connected vertex,
then $\varepsilon_G(v) = 1$. For example, all the vertices of a complete graph are
well-connected.

Recently, a number of topological indices involving vertex degree
and eccentricity were subject to a lot of mathematical as well as chemical studies.
A topological index of this type, introduced by Gupta et al. \cite{gup00}, was named as the connective eccentric index and was defined as $$C^\xi(G) = \sum\limits_{v \in V(G)}\frac{d_G(v)}{\varepsilon_G(v)}.$$ M. Ghorbani \cite{gho11} gave some bounds of connective eccentricity index and also computed this index for two infinite classes of dendrimers. De \cite{de12} reported some bounds for this index in terms of some graph invariants such as maximum and minimum degree, radius, diameter, first Zagreb index and first Zagreb eccentricity index, etc. In \cite{gho12a}, Ghorbani and Malekjani computed the eccentric connectivity index and the
connective eccentric index of an infinite family of fullerenes. In \cite{yu13}, Yu and
Feng also derived some upper or lower bounds for the connective eccentric index
and investigated the maximal and the minimal values of connective eccentricity
index among all $n$-vertex graphs with fixed number of pendent vertices.

In \cite{gho11}, Ghorbani showed that for a vertex transitive graph $G$, the
connective eccentric index is given by $$C^\xi(G) = \sum_{j = 1}^t |{A_j}|\frac{d_G(x_i)}{\varepsilon_G(x_i)},$$ where ${A_1},{A_2},\ldots,{A_t}$ are the orbits of $Aut(G)$ under its natural action on $V(G)$ and $x_i \in {A_i}, 1 \le i \le t$. In
particular, if $G$ is a regular graph, then ${C^\xi}(G) = \d\frac{n\delta}{r(G)}$, where $n$ is the number of vertices of $G$, which is a $\delta$-regular graph and $r(G)$ is the radius of $G$. Let, ${K_n},{C_n},{Q_m},{\Pi_m},{A_m}$ denote the complete graph with $n$
vertices, the cycle on $n$ vertices, $m$-dimensional hypercube,
$m$-sided prism and the $m$-sided antiprism respectively. It can be easily verified that the explicit formulae for the connective eccentric index of ${K_n},{C_n},{Q_m},{\Pi_m},{A_m}$ are as follows.

\begin{prop}
${C^\xi }({K_n}) = n(n - 1).$
\end{prop}

\begin{prop}
${C^\xi }({C_n}) = \frac{2n}{\lfloor n/2 \rfloor}.$
\end{prop}

\begin{prop}
${C^\xi }({Q_m}) = {2^m}.$
\end{prop}

\begin{prop}
${C^\xi }({\Pi _m}) = \left\{ \begin{array}{ll}
\frac{12m}{m + 2},&\mbox{when } m \mbox{ is even}\\[1mm]
\frac{12m}{m + 1},&\mbox{when } m \mbox{ is odd}.
\end{array}\right.$
\end{prop}

\begin{prop}
${C^\xi }({A_m}) = \left\{ \begin{array}{ll}
16,&\mbox{when } m \mbox{ is even}\\[1mm]
\frac{16m}{m + 1},&\mbox{when } m \mbox{ is odd}.
\end{array}\right.$
\end{prop}

Several studies on different topological indices related to graph operations of different
kinds are available in the literature \cite{ash11, gho12b, hos09, kha08, kha09}.

In this paper, first we calculate connective eccentric index of double graph and double cover
and hence the explicit formulae for the connective eccentric indices of join,
symmetric difference, disjunction and splice of graphs are obtained. For the
definitions and different results on graph operations, such as join, symmetric
difference, disjunction, etc., readers are referred to the book of Imrich and Klav\v{z}ar \cite{imr00}.

\section{Main Results}
In this section, first we define and then compute eccentric connectivity index of
double graph and double cover graph.

\subsection{Connective eccentric index of double graph and double cover}
Let us denote the double graph of a graph $G$ by $G^*$, which is constructed from two copies of
$G$ in the following manner \cite{alo86, hua12}. Let the vertex set of $G$ be $V(G) = \left\{
v_1, v_2, \ldots, v_n\right\}$, and the vertices of
${G^*}$ are given by the two sets $X = \left\{x_1, x_2, \ldots, x_n\right\}$ and $Y =
\left\{y_1, y_2, \ldots, y_n\right\}$. Thus for each vertex ${v_i} \in
V(G)$, there are two vertices ${x_i}$ and ${y_i}$ in $V(G^*)$. The double
graph $G^*$ includes the initial edge set of each copies of $G$, and for any edge
${v_i}{v_j} \in E(G)$, two more edges ${x_i}{y_j}$ and ${x_j}{y_i}$ are added.

\begin{figure}[h]
\begin{center}
\begin{tikzpicture}[place/.style={circle,draw=black!50,fill=black!20,
inner sep=0pt, minimum size = 5.2mm}]
{\footnotesize
\node[place] (1) {$v_1$};
\node[place] (2) [right=of 1] {$v_2$}
edge (1);
\node[place] (3) [right=of 2] {$v_3$}
edge (2);}
\end{tikzpicture}
\hspace{3cm}
\begin{tikzpicture}[place/.style={circle,draw=black!50,fill=black!20,
inner sep=0pt, minimum size = 5.2mm}]
{\footnotesize
\node[place] (1) {$x_1$};
\node[place] (2) [right=of 1] {$x_2$}
edge (1);
\node[place] (3) [right=of 2] {$x_3$}
edge (2);
\node[place] (4) [below=of 1] {$y_1$}
edge (2);
\node[place] (5) [below=of 2] {$y_2$}
edge (1) edge (3) edge (4);
\node[place] (6) [below=of 3] {$y_3$}
edge (2) edge (5);}
\end{tikzpicture}
\caption{The graph $P_3$ and its double graph $P_3^*$.}
\end{center}
\end{figure}
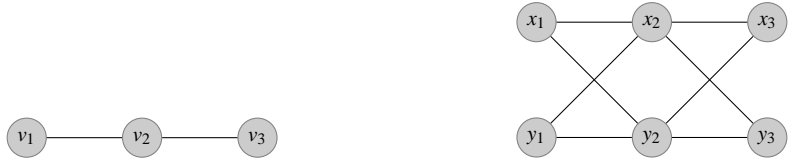

\begin{thm} The connective eccentric index of the double graph
${G^*}$ is given by ${C^\xi}(G^*) = 4{C^\xi}(G) - 2(n - 1){\left\|n - 1\right\|_G}$, where, ${\left\|n - 1\right\|_G}$ is the number of vertices with degree $n-1$, i.e., of eccentricity
one.
\end{thm}

\n\textit{Proof.} From the construction of double graph, it is clear that
$d_{G^*}(x_i) = d_{G^*}(y_i) = 2d_G(v_i)$, where $v_i \in V(G)$ and $x_i, y_i \in V(G^*)$ are the corresponding clone vertices of $v_i$. Also we can write $\varepsilon_{G^*}(x_i) = \varepsilon_{G^*}(y_i) = \varepsilon_G(v_i)$, when $\varepsilon_G(v_i) \ge 2$ and $\varepsilon_{G^*}(x_i) = \varepsilon_{G^*}(y_i) = \varepsilon_G(v_i) +
1$ when $\varepsilon_G(v_i) = 1$.

Thus the connective eccentric index of double graph $G^*$ is

\begin{eqnarray*}C^\xi(G^*)& = &\sum_{i = 1}^n\frac{d_{G^*}(x_i)}{\varepsilon_{G^*}(x_i)} + \sum_{i=1}^n \frac{d_{G^*}(y_i)}{\varepsilon_{G^*}(y_i)}\\
& = & 2\left[\sum_{\varepsilon_G(v_i) \ge 2}\frac{2d_G(v_i)}{\varepsilon_G(v_i)}+\sum_{\varepsilon_G(v_i)=1}\frac{2d_G(v_i)}{\varepsilon_G(v_i) + 1}\right]\\
&=&4\left[\sum_{\varepsilon_G(v_i)\ge 2}\frac{d_G(v_i)}{\varepsilon_G(v_i)}+\sum_{\varepsilon_G(v_i)=1}\frac{d_G(v_i)}{2}\right]\\
&=& 4\left[\sum_{\varepsilon_G(v_i)\ge 2}\frac{d_G(v_i)}{\varepsilon_G(v_i)}+ \sum_{\varepsilon_G(v_i)= 1}\frac{d_G(v_i)}{1}\right] - 2\sum_{\varepsilon_G(v_i)=1}d_G(v_i)\\
&=& 4C^\xi(G) - 2(n - 1){\left\|n - 1\right\|_G}.
\end{eqnarray*}
\qed

\bigskip
Let $G=(V,E)$ be a simple connected graph with $V = \left\{
{{v_1},{v_2},...,{v_n}} \right\}$. The extended double cover of $G$, denoted by
${G^{**}}$ is the bipartite graph with bipartition $(X,Y)$ where $X =
\left\{ {x_1},{x_2},...,{x_n}\right\}$ and $Y = \left\{
{y_1},{y_2},...,{y_n}\right\}$ in which ${x_i}$ and ${y_j}$ are adjacent if
and only if either ${v_i}$ and ${v_j}$ are adjacent in $G$ or $i = j$. For example,
the extended double cover of the complete graph is the complete bipartite graph.
This construction of the extended double cover was introduced by Alon \cite{alo86} in
1986.

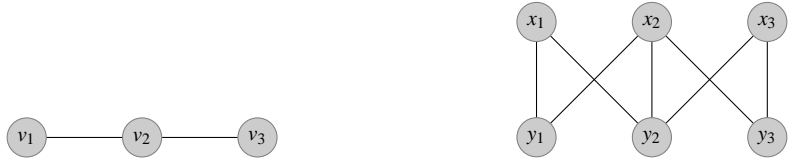
\begin{figure}[h]
\begin{center}
\begin{tikzpicture}[place/.style={circle,draw=black!50,fill=black!20,
inner sep=0pt, minimum size = 5.2mm}]
{\footnotesize
\node[place] (1) {$v_1$};
\node[place] (2) [right=of 1] {$v_2$}
edge (1);
\node[place] (3) [right=of 2] {$v_3$}
edge (2);}
\end{tikzpicture}
\hspace{3cm}
\begin{tikzpicture}[place/.style={circle,draw=black!50,fill=black!20,
inner sep=0pt, minimum size = 5.2mm}]
{\footnotesize
\node[place] (1) {$x_1$};
\node[place] (2) [right=of 1] {$x_2$};
\node[place] (3) [right=of 2] {$x_3$};
\node[place] (4) [below=of 1] {$y_1$}
edge (1) edge (2);
\node[place] (5) [below=of 2] {$y_2$}
edge (1) edge (2) edge (3);
\node[place] (6) [below=of 3] {$y_3$}
edge (2) edge (3);}
\end{tikzpicture}
\caption{The graph $P_3$ and its extended double cover $P_3^{**}$.}
\end{center}
\end{figure}

\begin{thm}
The connective eccentric index of the extended
double cover ${G^*}^*$ satisfies the inequality
\[
{C^\xi }({G^{**}}) \le 2{C^\xi }(G) - (n - 2){\left\| {n - 1} \right\|_G}.
\]
\end{thm}

\n\textit{Proof.} If $G$ is a graph with $n$ vertices and
$m$ edges then, from definition, the extended double cover graph ${G^{**}}$
consists of $2n$ vertices and $(n+2m)$ edges and $d_{G^{**}}(x_i) = d_{G^{**}}(y_i) = d_G(v_i) + 1$ and $\varepsilon _{G^*}(x_i) = \varepsilon _{G^*}(y_i) = \varepsilon_G(v_i) + 1$ for $i=1,2,\ldots,n$.

Thus the connective eccentric index of extended double cover graph ${G^{**}}$ is
given by
\begin{eqnarray*}
C^\xi(G^{**})& = &\sum_{i = 1}^n
\frac{d_{G^{**}}(x_i)}{\varepsilon_{G^{**}}(x_i)} +
\sum_{i = 1}^n\frac{d_{G^*}(y_i)}{\varepsilon_{G^*}(y_i)}\\
&=& 2\left[\sum_{\varepsilon_G(v_i) \ge 2}
\frac{d_G(v_i) + 1}{\varepsilon_G(v_i) + 1} +
\sum_{\varepsilon_G(v_i) = 1}\frac{d_G(v_i) + 1}{2}
\right]\\
&\le & 2\sum_{\varepsilon_G(v_i) \ge 1}\frac{d_G(v_i)}{\varepsilon_G(v_i)} + \sum_{\varepsilon_G(v_i) = 1}\left\{d_G(v_i) + 1\right\}\\
&= & 2\left[\sum_{\varepsilon_G(v_i)\ge 1}\frac{d_G(v_i)}{\varepsilon_G(v_i)} +  \sum_{\varepsilon_G(v_i) = 1}d_G(v_i)\right] - \sum_{\varepsilon_G(v_i) = 1}d_G(v_i) +\sum_{\varepsilon_G(v_i) = 1} 1\\
& = & 2C^\xi(G) - (n - 1){\left\|n - 1\right\|_G} + {\left\|n - 1\right\|_G}\\
&=& 2C^\xi(G) - (n - 2){\left\|n - 1\right\|_G}.
\end{eqnarray*}
\qed

\bigskip

Now some exact formula for the eccentric connectivity index of joins, symmetric
difference, disjunction, splice graphs are presented.

\subsection{Join}
The join $G_1 + G_2$ of two graphs $G_1$ and $G_2$ with disjoint vertex
sets $V_1, V_2$ and edge sets $E_1, E_2$, is the graph union $G_1 \cup G_2$ together with all the edges joining ${V_1}$ and ${V_2}$, i.e., $G_1 +G_2$ consists of the
vertex set $V_1\cup V_2$ and edge set $E_1\cup E_2 \cup \{ xy:x \in V_1, y \in V_2\}$.

\begin{thm} Let ${G_1}$ and ${G_2}$ be two graphs without
well-connected vertices. Then,
$
{C^\xi }({G_1} + {G_2}) = |E({G_1})| + |E({G_2})| + |V({G_1})||V({G_2})|.
$
\end{thm}

\n\textit{Proof.} Let $n_1$ and $n_2$ be the numbers of vertices in $G_1$ and $G_2$ respectively. Thus, $|V({G_1} + {G_2})| = |{V_1}| + |{V_2}| = {n_1} + {n_2}$. For vertices $u \in {V_1}$ and $v \in {V_2}$, it holds that $d_{G_1+G_2}(u) = d_{G_1}(u) + n_2$ and $d_{G_1+G_2}(v) =
d_{G_2}(v) + n_1$. Since none of ${G_1}$ and ${G_2}$ contains
well connected vertices, then for every $u \in V({G_1} + {G_2}), \varepsilon_{G_1+G_2}(u) = 2$. So, from definition of connective eccentric index, we have
\begin{eqnarray*}
C^\xi(G_1+G_2)& = &\sum_{v\in V(G_1)}\frac{d_{G_1+G_2}(v)}{\varepsilon_{G_1+G_2}(v)}+ \sum_{v\in
V(G_2)}\frac{d_{G_1+G_2}(v)}{\varepsilon_{G_1+G_2}(v)}\\
&= &\frac{1}{2}\left[\sum_{v \in V(G_1)}\left\{d_{G_1}(v)+n_2 \right\}+\sum_{v\in V(G_2)}\left\{d_{G_2}(v)+n_1\right\}\right]
\end{eqnarray*}
from where the desired result follows.
\qed

\bigskip
Now we generalize the above result for $n$ disjoint graphs $G_1,G_2,\ldots,{G_n}$.

\begin{thm} Let ${G_1},{G_2},\ldots,{G_n}$ be $n$ graphs with disjoint vertex sets
${V_i} = V({G_i})$ and edge sets ${E_i} = E({G_i}), 1 \le i \le n$, without
well-connected vertices. Then,
\[
{C^\xi }({G_1} + {G_2} + ... + {G_n}) = \sum\limits_{i = 1}^n {|E({G_i})|}  +
\frac{1}{2}{\left\{ {\sum\limits_{i = 1}^n {|V({G_i})|} } \right\}^2} -
\frac{1}{2}\sum\limits_{i = 1}^n {|V({G_i}){|^2}}.
\]

\end{thm}

\n\textit{Proof.} From definition of join, we have
$
d_{G_1 + G_2 + ... + G_n}(v_j) = d_{G_i}(v_j) + |V(G_1 +G_2 + ... + G_n)| - |V(G_i)|$ and since none of $G_1, G_2, ..., G_n$ contains any well connected
vertex, we have $\varepsilon_{G_1 + G_2 + ... + G_n}(v_j) = 2$ for all $v_j\in V(G_i)$.

Thus the connective eccentric index of $G_1+ G_2+ ... + G_n$ is given by
\begin{eqnarray*}
C^\xi(G_1 + G_2 + \ldots + G_n) &=& \sum_{i = 1}^n\sum_{j
= 1}^{|V({G_i})|}\frac{d_{G_1 + G_2 + ... +G_n}(v_j)}{\varepsilon_{G_1 + G_2 + ... +G_n}(v_j)}\\
&=& \sum_{i = 1}^n\sum_{j = 1}^{|V(G_i)|}\frac{d_{G_i}(v_j) + \left|V(G_1 + G_2 + ... + G_n)\right| - \left|V(G_i)\right|}{2} \\
& = & \frac{1}{2}\bigg[|E(G_1)| + |E(G_2)| + ... +|E(G_n)|\\
&&+ |V(G_1+ G_2+... + G_n)|\left\{|V(G_1)| + |V(G_2)| + ... \right.\\
&&\left. + |V(G_n)|\right\}\left. - \left\{ |V(G_1)|^2 + V(G_2)|^2 +... + |V(G_n)|^2\right\}
\right]
\end{eqnarray*}
from where the desired result follows.
\qed

\bigskip
The following corollaries are direct consequences of the theorem.

\begin{cor}
If $nG$ denotes the join of $n$ copies of $G$, then
\[
{C^\xi }(nG) = n|E(G)| + \frac{{n(n - 1)}}{2}|V(G){|^2}.
\]
\end{cor}

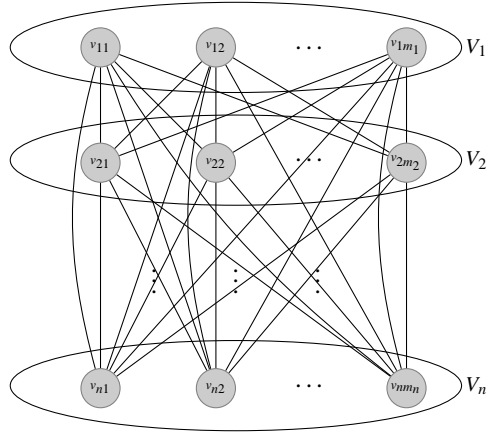
\begin{figure}[h]
\begin{center}
\begin{tikzpicture}[place/.style={circle,draw=black!50,fill=black!20,
inner sep=0pt, minimum size = 5.2mm}, bend angle=15]
\draw (1.8,0) ellipse (30mm and 6mm);
\draw (1.8,-1.5) ellipse (30mm and 6mm);
\draw (1.8,-4.5) ellipse (30mm and 6mm);
\draw (2.8,0) node{\ldots};
\draw (2.8,-1.5) node{\ldots};
\draw (1.8,-3) node{\vdots \hspace{1cm}\vdots \hspace{1cm}\vdots};
\draw (2.8,-4.5) node{\ldots};
\draw (5,0) node{\footnotesize $V_1$};
\draw (5,-1.5) node{\footnotesize $V_2$};
\draw (5,-4.5) node{\footnotesize $V_n$};
{\tiny
\node[place] (1) {$v_{11}$};
\node[place] (2) [right=of 1] {$v_{12}$};
\node[place] (3) [right=of 2, xshift=1cm] {$v_{1m_1}$};
\node[place] (4) [below=of 1] {$v_{21}$}
edge (1) edge (2) edge (3);
\node[place] (5) [below=of 2] {$v_{22}$}
edge (1) edge (2) edge (3);
\node[place] (6) [below=of 3] {$v_{2m_2}$}
edge (1) edge (2) edge (3);
\node[place] (7) [below=of 1, yshift=-3cm]{$v_{n1}$}
edge[bend left] (1) edge (2) edge (3) edge (4) edge (5) edge (6);
\node[place] (8) [below=of 2, yshift=-3cm] {$v_{n2}$}
edge (1) edge[bend left] (2) edge (3) edge (4) edge (5) edge (6);
\node[place] (9) [below=of 3, yshift=-3cm] {$v_{nm_n}$}
edge[bend angle=10, bend left] (1) edge (2) edge[bend left] (3) edge (4) edge (5) edge (6);
}
\end{tikzpicture}
\caption{The complete $n$-partite graph.}
\end{center}
\end{figure}

\begin{cor} Let $G = K_{m_1, m_2, ..., m_n}$ be the complete
$n$-partite graph having $v = |V(G)|$ number of vertices. Here the vertex
set can be partitioned into $n$ subsets ${V_1},{V_2},...,{V_n}$ such that ${V_1}
\cup {V_2} \cup ... \cup {V_n} = V(G)$ and ${V_i} \cap {V_j} = \phi,i \ne j, 1 \le i,j \le n$. From definition of join it is clear
that $G = K_{m_1,m_2,...,m_n}$ is the join of $n$ empty graphs with
${m_1},{m_2},...,{m_n}$  number of vertices. Then the
connective eccentric index of $G = K_{m_1, m_2,..., m_n}$ is given by
\[
{C^\xi }(K_{m_1, m_2, ..., m_n}) = \sum_{{1 \le i\le n\atop 1\le j \le n}\atop i \ne j}
{{m_i}{m_j}}.
\]
\end{cor}

\subsection{Symmetric difference}
Let ${G_1}$ and ${G_2}$ be two graphs with vertex sets $V({G_1})$ and $V({G_2})$ and
edge sets $E({G_1})$ and $E({G_2})$. Then the symmetric difference of
${G_1}$ and ${G_2}$, denoted by ${G_1} \oplus {G_2}$ is the graph with vertex set
$V({G_1}) \times V({G_2})$ in which any two vertices $({u_1},{u_2})$ is adjacent
to $({v_1},{v_2})$ whenever ${u_1}$ is adjacent to ${v_1}$ in $G_1$ or
${u_2}$ is adjacent to ${v_2}$ in $G_2$, but not both. From definition
of symmetric difference, the degree of a vertex $({u_1},{u_2})$ of ${G_1} \oplus
{G_2}$ is given by \cite{imr00}
\[
d_{G_1\oplus G_2}({u_1},{u_2}) = |V({G_2})|{d_{{G_1}}}({u_1}) +
|V({G_1})|{d_{{G_2}}}({u_2}) - 2{d_{{G_1}}}({u_1}){d_{{G_2}}}({u_2}).
\]

\begin{thm} The connective eccentric index of the symmetric
difference ${G_1} \oplus {G_2}$ of two graphs ${G_1}$ and ${G_2}$ is given by
\[
{C^\xi }({G_1} \oplus {G_2}) = |E({G_1})||V({G_2}){|^2} +
|E({G_2})||V({G_1}){|^2} - 4|E({G_1})||E({G_2})|
\]
where, none of ${G_1}$ and ${G_2}$ contains well-connected vertices.
\end{thm}

\n\textit{Proof.} Since the distance between any two vertices of a symmetric
difference cannot exceed two, if none of the components contains well-connected
vertices, the eccentricity of all vertices is constant and equal to two, i.e., ${\varepsilon _{{G_1} \oplus {G_2}}}({u_1},{u_2}) = 2$, for all
vertices $({u_1},{u_2})$ \cite{ash11, gho12b}.

Thus the connective eccentric index of symmetric difference ${G_1}
\oplus {G_2}$ of two graphs ${G_1}$ and ${G_2}$ is given by
\begin{eqnarray*}
C^\xi(G_1\oplus G_2)&=& \sum_{(u_1,u_2) \in V(G_1\oplus G_2)}\frac{d_{G_1\oplus G_2}(u_1,u_2)}{\varepsilon_{G_1\oplus G_2}(u_1,u_2)}\\
&=& \frac{1}{2}\sum_{u_1\in V(G_1)}\sum_{u_2\in V(G_2)}\left[ |V(G_2)|d_{G_1}(u_1) + |V(G_1)|d_{G_2}(u_2) -2d_{G_1}(u_1)d_{G_2}(u_2) \right]\\
&=& \frac{1}{2}\left[ |V(G_2)|^2\sum_{u_1 \in V(G_1)}d_{G_1}(u_1) + |V(G_1)|^2\sum_{u_2\in V(G_2)}d_{G_2}(u_2)\right.\\
&&\left.  - 2\sum_{u_1\in V(G_1)}d_{G_1}(u_1)\sum_{u_2\in V(G_2)}d_{G_2}(u_2)\right]
\end{eqnarray*}
from where the desired result follows.
\qed

\subsection{Disjunction}
The disjunction ${G_1} \vee {G_2}$of two graphs ${G_1}$and ${G_2}$is the graph
with vertex set $V({G_1}) \times V({G_2})$ in which $({u_1},{u_2})$ is adjacent to
 $({v_1},{v_2})$ whenever ${u_1}$ is adjacent with ${v_1}$ in $G_1$
or ${u_2}$ is adjacent with ${v_2}$ in $G_2$. Obviously, the degree of a
vertex $({u_1},{u_2})$ of ${G_1} \vee {G_2}$ is given by \cite{ash11, gho12b}
\[
{d_{({G_1} \vee {G_2})}}({u_1},{u_2}) = \left[ {|V({G_2})|{d_{{G_1}}}({u_1}) +
|V({G_1})|{d_{{G_2}}}({v_2}) - {d_{{G_1}}}({u_1}){d_{{G_2}}}({v_2})} \right].
\]

\begin{thm}The connective eccentric index of the disjunction ${G_1}
\vee {G_2}$ of two graphs${G_1}$ and ${G_2}$ is given by
\[
{C^\xi }({G_1} \vee {G_2}) = |E({G_1})||V({G_2}){|^2} + |E({G_2})||V({G_1}){|^2}
- 2|E({G_1})||E({G_2})|
\]
where none of ${G_1}$ and ${G_2}$ contains well-connected vertices.
\end{thm}

\n\textit{Proof.} Since the distance between any two vertices of a disjunction
cannot exceed two, if none of the components contains well-connected vertices,
the eccentricity of all vertices is constant and equal to two \cite{ash11, gho12b}. Then the
connective eccentric index of the disjunction ${G_1} \vee {G_2}$ of two
graphs ${G_1}$ and ${G_2}$ is computed as
\begin{eqnarray*}
C^\xi(G_1\vee G_2) &=& \sum_{(u_i,v_j) \in V(G_1\vee G_2)}\frac{d_{G_1\vee G_2}(u_1,u_2)}{\varepsilon_{G_1\vee G_2}(u_1,u_2)}\\
&=& \frac{1}{2}\sum\limits_{{u_1} \in V({G_1})} {\sum\limits_{{v_2} \in V({G_2})}
{\left[ {|V({G_2})|{d_{{G_1}}}({u_1}) + |V({G_1})|{d_{{G_2}}}({u_2}) -
{d_{{G_1}}}({u_1}){d_{{G_2}}}({u_2})} \right]} }
\end{eqnarray*}
from where the desired result follows.
\qed

\subsection{Splice}
Let ${G_1} = ({V_1},{E_1})$ and ${G_2} = ({V_2},{E_2})$ are two graphs ${V_1} \cap
{V_2} = \phi $. Let ${v_1} \in {V_1}$ and ${v_2} \in {V_2}$ be two given vertices
of ${G_1}$ and ${G_1}$ respectively. A splice of ${G_1}$ and ${G_1}$ at the
vertices ${v_1}$ and ${v_2}$ is obtained by identifying the vertices ${v_1}$
and ${v_2}$ in the union of ${G_1}$ and ${G_1}$ and is denoted
by $S({G_1},{G_2},{v_1},{v_2})$ \cite{dos05}. Different topological indices of the splice
graphs have already been computed \cite{sha13}. Let $u \in V(S)$, then ${\varepsilon
_1}(u)$ denotes the eccentricity of $u$ as a vertex of ${G_1}, {\varepsilon
_2}(u)$ denotes the eccentricity of $u$ as a vertex of ${G_2}$ and
$\varepsilon (u)$ denotes the eccentricity of $u$ as a vertex of $G$.
 Let ${d_S}(x)$ be the degree of the vertex $x \in V(S)$. Then the connective
eccentric index of $S({G_1},{G_2},{v_1},{v_2})$ is computed as follows.

\begin{thm} The connective eccentric index of splice of  ${G_1}$
and ${G_2}$ is given by
\[
{C^\xi }(S) = \sum\limits_{x \in {V_1}} {\frac{{{d_S}(x)}}{{max\left\{
{{d_{{G_1}}}(x,{u_1}) + \varepsilon ({u_2})\varepsilon (x)} \right\}}} + }
\sum\limits_{y \in {V_2}} {\frac{{{d_S}(y)}}{{max\left\{ {{d_{{G_2}}}(y,{u_2}) +
{\varepsilon _1}({u_1})\varepsilon (y)} \right\}}}}.
\]
\end{thm}

\n\textit{Proof.} For any vertex $x \in {V_1}$, Sharafdini and Gutman \cite{sha13} showed
that
\[
\varepsilon (x) = \max\left\{ {{d_{{G_1}}}(x,{u_1}) + {\varepsilon
_2}({u_2}){\varepsilon _1}(x)} \right\}.
\]
Similarly, for any vertex $y \in {V_2}$,
\[
\varepsilon (y) = \max\left\{ {{d_{{G_2}}}(y,{u_2}) + {\varepsilon _1}({u_1})\varepsilon (y)} \right\}.
\]
Thus, the desired result follows from the definition of connective eccentric index.
\qed

\end{document}